\documentclass{amsart}

\usepackage{amsfonts}
\usepackage[english]{babel}
\def\g{\gamma}
\def\G{\Gamma}
\def\d{\delta}
\def\a{\alpha}
\def\b{\beta}
\def\p{\varphi}
\def\e{\varepsilon}
\def\l{\lambda}

\def\s{\sigma}
\def\k{\kappa}
\def\o{\omega}
\def\O{\Omega}
\def\z{\zeta}
\def\r{\rho}

\def\W{\mathcal W}
\def\SH{\mathcal S_H}
\def\D{\mathcal D}

\def\R{{\mathbb R}}
\def\C{{\mathbb C}}
\def\N{{\mathbb N}}
\def\Z{{\mathbb Z}}
\def\Re{\mbox{Re }}
\def\Im{\mbox{Im }}

\DeclareMathOperator{\spec}{sp}
\DeclareMathOperator{\Res}{Res}

\DeclareMathOperator{\supp}{supp }
\DeclareMathOperator{\dist}{dist}

\newtheorem{Th}{Theorem}
\newtheorem{Pro}{Proposition}

\begin{document}

\title{Analogues of Fourier quasicrystals for a strip}

\author{Sergii Yu.Favorov}

\address{Sergii Favorov,
\newline\hphantom{iii}  V.N.Karazin Kharkiv National University
\newline\hphantom{iii} Svobody sq., 4, Kharkiv, Ukraine 61022}
\email{sfavorov@gmail.com}

\maketitle {\small
\begin{quote}
\noindent{\bf Abstract.}
We study a certain family of discrete measures with unit masses on a horizontal strip as an analogue of Fourier quasicrystals on the real line. We prove a one-to-one correspondence
between supports of  measures from this family and zero sets of exponential polynomials with imaginary frequencies. This result is the special case of a general result on measures
whose supports correspond to zero sets of absolutely convergent Dirichlet series with bounded spectrum.

\medskip


\end{quote}
}

\medskip

   \section{Introduction}\label{S1}
   \bigskip

 A measure $\mu$ with locally finite support (that is, its intersection with any compact set is finite) is called a Fourier quasicrystal if $\mu$ is a temperate distribution, its Fourier transform
in the sense of distributions $\hat\mu$ is also a measure with locally finite support, and both measures $|\mu|$, $|\hat\mu|$ are temperate distributions.
Here and below $|\nu|(E)$ means the variation of the complex measure $\nu$ on the set $E$.

 In \cite{OU}, \cite{OU1} A.Olevskii and A.Ulanovskii proved that a measure $\mu_A$  of the form
\begin{equation}\label{a}
\mu_A=\sum_n \d_{a_n},\quad A=\{a_n\}\subset\R,
\end{equation}
  is a Fourier quasicrystal  if and only if its support $A=\{a_n\}$ is the zero set of an exponential polynomial
\begin{equation}\label{s}
P(z)=\sum_{1\le j\le N} q_j  e^{2\pi i\o_j z}, \qquad q_j\in\C,\quad \o_j\in\R.
\end{equation}
Here $\d_x$ means the unit mass at the point $x$.

From now on, we always assume, that each point can occur in $A$ a finite number of times, so in fact $A$ is a multiset.

Fourier quasicrystals are currently being studied very actively. Many works are devoted to the study of the properties of Fourier quasicrystals (see, for example, the collections of papers \cite{D}, \cite{Q},
and the basic paper \cite{L1}). Fourier quasicrystals find applications in modern physics, where they serve as mathematical models of certain atomic structures. A nontrivial example of
a Fourier quasicrystal of the form \eqref{a}, whose support has only a finite intersection with any arithmetical progression, was found by P.Kurasov and P.Sarnak \cite{KS}.
Let us mention also the paper \cite{G}, which actually considers pairs $(\mu,\hat\mu)$.
\medskip

In \cite{F2}, \cite{F3} we extended Olevskii--Ulanovskii's scheme to zero sets $A=\{a_n\}\subset\R$ of absolutely convergent Dirichlet series with only real zeros of the form
\begin{equation}\label{S}
Q(z)=\sum_{\o\in\O} q_\o  e^{2\pi i\o z}, \qquad q_\o\in\C,\quad\sum_{\o\in\O} |q_\o|<\infty,\quad \O\subset\R\,\,\text{is bounded}.
\end{equation}
 Such sums are natural extensions  of the finite sums \eqref{s}. We proved that the Fourier transform of the measure $\mu_A$ defined in \eqref{a} is also a measure of the form
\begin{equation}\label{b}
  \hat\mu_A=\sum_{\g\in\G}b_\g\d_\g,\quad \G\subset\R\ \text{ is countable}.
\end{equation}
 Besides, we  found  necessary and sufficient conditions on
$\hat\mu_A$ for the set $A$ to be the zero set of Dirichlet series \eqref{S}.

The  zero sets $A=\{a_n\}$ of  exponential polynomials \eqref{s} and Dirichlet series \eqref{S}
with zeros $a_n$ either on the real line or on a horizontal strip of finite width are almost periodic sets in the sense of M.Krein and B.Levin  \cite[App.VI]{L}.
In particular, $A$ is {\sl translation bounded}, i.e., the quantities\footnote{$\# A$ means a number of points of  the finite set $A$; if $A$ is a multiset, points are counted according to their multiplicities.}
$$
\#\{n:\,t<\Re a_n<t+1\}
$$
are bounded uniformly in $t\in\R$.
On the other hand, the measure $|\hat\mu|$ on $\R$ is temperate if and only if the variation $|\hat\mu|(-r,r)$ grows polynomially as $r\to\infty$ (see \cite{F1}).
 Having this in mind, we can reformulate the above result of Olevskii and Ulanovskii as follows:
\medskip

{\it A locally finite set $A=\{a_n\}\subset\R$ is the zero set of an exponential polynomial \eqref{s} if and only if $A$ is almost periodic, $\hat\mu_A$ has the form \eqref{b}
with locally finite support $\G=\{\g\in\R:\,b_\g\neq0\}$, and
$$
\log\sum_{|\g|<r}|b_\g|=O(\log r),\qquad\ r\to\infty.
$$}
\medskip

Here we  apply the method of Olevskii--Ulanovskiie to the measures of the form \eqref{a} where  we replace the condition $A\subset\R$ with $A\subset\SH:=\{z=x+iy:\,|y|\le H\}$.
  We slightly change the definition of $\hat\mu_A$ and prove that the above results on zeros of exponential sums and Dirichlet series are valid in this case as well.

Namely, let $\mu$ be a measure with support in $\SH$. Denote
$$
  M_\mu(r):=|\mu|(\{z\in\SH:\,|z|\le r\}),
$$
and suppose that
\begin{equation}\label{log}
    \log M_\mu(r)=O(\log r)\qquad r\to\infty.
\end{equation}
Let $\D$ be the space of all $C^\infty$-functions with compact support. For $\p\in\D$ put
$$
   \hat\p^c(z):=\int_{\R}\phi(t)e^{-2\pi i zt}dt,
$$
and define {\sl the c-Fourier transform} of a measure $\mu$  by the equality
\begin{equation}\label{h}
(\hat\mu^c,\p)=(\mu,\hat\p^c)=\int\hat\p^c(z)\mu(dz).
\end{equation}
The right-hand side is finite  (see Section \ref{S2}). In the general case $\hat\mu^c$ is an element of the space $\D'$.

\begin{Th}\label{T1}
Let $A=\{a_n\}\subset\SH$  be the zero set of some Dirichlet series \eqref{S} and
\begin{equation}\label{a0}
\mu_A=\sum_n\d_{a_n}.
\end{equation}
 Then $A$ is an almost periodic set and the c-Fourier transform $\hat\mu^c_A$ is a pure point measure \eqref{b} with coefficients $b_\g$ satisfying the conditions
\begin{equation}\label{var}
\log\sum_{|\g|<r}|b_\g|=O(r),\qquad r\to\infty,
\end{equation}
and
\begin{equation}\label{neig}
\sum_{0<|\g|<1}|b_\g/\g|<\infty.
\end{equation}
If $A$ is the zero set of an exponential polynomial \eqref{s}, then $\G$ in \eqref{b} is locally finite.
\end{Th}

\begin{Th}\label{T2}
Let $A=\{a_n\}\subset\SH$ be an almost periodic set,  let $\mu_A$ be a measure defined in \eqref{a0}, and let the c-Fourier transform $\hat\mu^c_A$ be a pure point measure \eqref{b} with coefficients $b_\g$
satisfying  \eqref{var}  and \eqref{neig}. Then there exists a Dirichlet series \eqref{S} with  zero set $A$. If, in addition, $\G$ is locally finite,
then the Dirichlet series is  an exponential polynomial \eqref{s}.
\end{Th}

The article is structured as follows.

In sections \ref{S2} and \ref{S3} we present the properties of the c-Fourier transform and  Dirichlet series we need, respectively.

In section \ref{S4} we give definitions of almost periodic measures and sets and describe some of their properties.

In Section \ref{S5} we consider  entire functions of exponential growth with an almost periodic zero set
 and prove some formulas for them and their logarithmic derivative.

In Sections \ref{S6} and \ref{S7} we prove  Theorem \ref{T1} and \ref{T2}, respectively.

\bigskip
\section{The Fourier transform on the complex plane}\label{S2}
\bigskip

 For $\psi\in L^1(\R)$ put
$$
   \hat\psi(x)=\int_\R \psi(t)e^{-2\pi ixt}dt.
$$
 If
\begin{equation}\label{e}
\psi(t)=O\left(e^{-K|t|}\right),\qquad K>2\pi H,\qquad |t|\to\infty,
\end{equation}
 then the function
\begin{equation}\label{ext}
   \hat\psi^c(z)=\int_{\R}\psi(t)e^{-2\pi i zt}dt=\widehat{(\psi(t)e^{2\pi yt})}(x),\quad z=x+iy,
\end{equation}
 is the holomorphic extension of $\hat\psi(x)$ to a neighborhood of $\SH$. Extending to the strip the well-known equality for real $x$, we get for every $\p\in\D$ and $\psi\in L^1(\R)$ which satisfies \eqref{e}
 \begin{equation}\label{com}
\widehat{(\p\star\psi)^c}(z)=\hat\p^c(z)\hat\psi^c(z),\quad \forall z\in\SH.
\end{equation}
If $\psi$ is an $m$-continuously differentiable function such that \eqref{e} holds for $\psi$ and all its derivatives $\psi^{(k)},\,k\le m$, then
\begin{equation}\label{phi}
  \hat\psi^c(z)=(2\pi iz)^{-m}\int_\R\psi^{(m)}(t)e^{-2\pi itz}dt, \quad z\in\SH\setminus\{0\}.
\end{equation}
For all $z\in\SH$ the integral is finite, hence,
\begin{equation}\label{N}
  |\hat\psi^c(z)|\le C(\max\{1,|z|\})^{-m},\qquad C=C(\psi,m)<\infty.
\end{equation}
Applying \eqref{phi} to $t^k\psi(t)$, we get  inequalities similar to \eqref{N} for  $(\hat\psi^c)^{(k)}(z)$ for all $k\in\N$.

Note that \eqref{log} is equivalent to the estimate
\begin{equation}\label{lo}
M_\mu(r)=O(r^N), \qquad r\to\infty.
\end{equation}
with some $N\in\N$. Using \eqref{N} for $\p\in\D,\,m=N+1$ and integrating by parts, we get
$$
 \left|\int_{z\in\SH,|z|>R}\hat\p^c(z)\mu(dz)\right|\le C\int_{z\in\SH,|z|>R}|z|^{-N-1}|\mu|(dz)\le C(N+1)\int_R^\infty r^{-N-2}M_\mu(r)dr.
$$
The last integral tends to zero as $R\to\infty$. Therefore the integral in \eqref{h} is finite, and the distribution $\hat\mu^c$ is well-define.

\begin{Pro}\label{P1}
Let  $\mu$ be a measure on $\SH$ satisfying \eqref{lo}, and let $\psi$ be a $C^\infty$-function such that all its derivatives satisfy \eqref{e} with $K>2\pi H$.
If  $\hat\mu^c$ is a measure on $\R$ such that
$$
  |\hat\mu^c|(-r,r)=O(e^{Lr}),\qquad 2\pi H<L<K,\qquad r\to\infty,
$$
 then the equality \eqref{h} remains valid with the function $\psi$ in place of $\p\in\D$..
\end{Pro}

{\bf Proof}.  Let $\eta(t)\in\D$ such that $\eta(0)=1$. Set $\eta_n(t)=\eta(t/n)$.
Replacing $\psi$ with $\eta_n\psi$ in \eqref{ext} and applying Lebesgue's Dominate Convergence Theorem, we get for all $z\in\SH$
$$
(\widehat{\eta_n\psi})^c(z)\to \hat\psi^c(z),\qquad n\to\infty.
$$
Next, it is easy to see that
$$
|(\eta_n(t)\psi(t))^{(N+1)}|\le Ce^{-K|t|}, \quad C<\infty,
$$
uniformly in $n$. Hence, \eqref{N} with $m=N+1$  and with $\eta_n\psi$ in place of $\psi$ imply
$$
  |(\widehat{\eta_n\psi})^c(z)|\le C'(\max\{1,|z|\})^{-N-1},\quad z\in\SH,\quad C'=C'(\eta,\psi,N)<\infty.
$$
Since the integral
$$
 \int|(\widehat{\eta_n\psi})^c(z)||\mu|(dz)\le C'\int_0^\infty(\max\{1,r\})^{-N-1}M_\mu(dr)=C'(N+1)\int_1^\infty r^{-N-2}M_\mu(r)dr
$$
is finite, we can apply  Lebesgue's Dominate Convergence Theorem and get
\begin{equation}\label{ps1}
\lim_{n\to\infty}(\mu,(\widehat{\eta_n\psi})^c)=\lim_{n\to\infty}\int(\widehat{\eta_n\psi})^c(z)\mu(dz)=\int\hat\psi^c(z)\mu(dz)=(\mu,\hat\psi^c).
\end{equation}
On the other hand, taking into account \eqref{e} and integrating by parts, we get
$$
 \int_\R|\eta_n(t)\psi(t)||\hat\mu^c|(dt)\le C\int_\R e^{-K|t|}|\hat\mu^c|(dt)\le KC\int_o^\infty e^{-Kr}|\hat\mu^c|(-r,r)dr<\infty.
$$
Applying  Lebesgue's Dominate Convergence once more, we obtain
\begin{equation}\label{ps2}
\lim_{n\to\infty}(\hat\mu^c,\eta_n\psi)=\lim_{n\to\infty}\int_\R\eta_n(t)\psi(t)\hat\mu^c(dt)=\int_\R\psi(t)\hat\mu^c(dt)=(\hat\mu^c,\psi).
\end{equation}
Since $\eta_n\psi\in\D$, we see that \eqref{ps1}, \eqref{ps2} and the definition of $\hat\mu^c$  imply the statement of the Proposition.$\square$

\bigskip
\section{Dirichlet series on the real line}\label{S3}
\bigskip

 Denote by $\W$ the algebra of all  Dirichlet series on $\R$
$$
Q(x)=\sum_n q_ne^{2\pi i\o_n x},\,\o_n\in\R,
$$
with the  finite Wiener's norm $\|Q\|_\W=\sum_n|q_n|$. The spectrum of $Q\in\W$ is the set $\spec Q=\{\o_n:\,q_n\neq0\}$.
For each $Q\in\W$ and each analytic function $h(z)$ on a neighborhood of the set $\overline{\{Q(x):\,x\in\R\}}$ we have $h(Q(x))\in\W$ (see \cite{R1}, Ch.VI). In particular, $\exp Q\in\W$,
if $\inf_\R|Q(x)|>0$ then $1/Q\in\W$, and if $\sup_\R|Q(x)|<1$ then $\log(1+Q(x))\in\W$.

When the spectrum $\spec Q$ is non-negative and locally finite $\spec Q$ the set $\spec Q^n\cap(0,r)$ for all $n,\, r<\infty$ is a subset of all sums of
at most $\frac{r}{\min\{\spec Q\setminus\{0\}\}}$ possibly identical elements of $\spec Q\cap(0,r)$.
Therefore, the sets $\spec Q^n\cap(0,r)$ for all $n$ are subsets of one finite set, which depends only on $r$, and the equalities
$$
\exp{Q(x)}=\sum_{n=0}^\infty Q^n(x)/n!\quad\text{for}\quad \|Q\|_\W<\infty,\qquad \log(1+Q(x))=\sum_{n=1}^\infty (-1)^{n-1} Q^n(x)/n\quad\text{for}\quad \|Q\|_\W<1,
$$
imply the following proposition:
\begin{Pro}\label{P2}
If the spectrum of $Q\in\W$  is non-negative and  locally finite, then  so is the spectrum $\spec\exp Q$. If, in addition,  $\|Q\|_\W<1$   and $0\not\in\spec Q$,
then $\spec\log(1+Q)$  is strictly positive and locally finite too.
\end{Pro}

\bigskip
\section{Almost periodic  measures and sets on a strip}\label{S4}
\bigskip

{\bf Definition 1} (see \cite{B}).
 A continuous function $g(x)$ on the real line $\R$
is  \emph{almost periodic}  if for any  $\e>0$ the set of its $\e$-almost periods
  $$
E_\e= \{\tau\in\R:\,\sup_{x\in\R}|g(x+\tau)-g(x)|<\e\}
  $$
is relatively dense, i.e., $E_\e\cap(t,t+L)\neq\emptyset$ for all $t\in\R$ and some $L$ depending on $\e$.

For example, any function $Q\in\W$ is almost periodic.
\medskip

{\bf Definition 2} (see \cite{B}, App.II).
 A continuous function $g(z)$ on the open strip
$$
S_{(\a,\b)}=\{z=x+iy\in\C:\,-\infty\le\a<y<\b\le+\infty\}
$$
is  \emph{almost periodic}  if for any $\a',\b'$ such that $\a<\a'<\b'<\b$ and
  $\e>0$ the set of $\e$-almost periods
  $$
E_{\a',\b',\e}= \{\tau\in\R:\,\sup_{x\in\R,\a'\le y\le\b'}|g(x+\tau+iy)-g(x+iy)|<\e\}
  $$
is relatively dense, i.e., $E_{\a',\b',\e}\cap(t,t+L)\neq\emptyset$ for all $t\in\R$ and some $L$ depending on $\a',\b',\e$.
\medskip

Every holomorphic function, which is bounded on $S_{(\a,\b)}$ and almost periodic on one line $\Im z=c,\,\a<c<\b$ is almost periodic in  $S_{(\a,\b)}$ (see \cite[Part 2, Ch.1]{Le}).
For example, every Dirichlet series with finite Wiener's norm and bounded spectrum \eqref{S} extends as an entire almost periodic function on the plane $\C=S_{(-\infty,+\infty)}$.
\medskip

{\bf Definition 3} (\cite{R}).
 A  measure $\mu$   is \emph{almost periodic}  on the open strip $S_{(\a,\b)}$  if for  any continuous function $\p(z)$ with compact support
  contained in the set $\{z\in\C:\,\a-\a'<y<\b-\b'\}$,  the convolution  $\p\star\mu(z)=\int\p(z-w)\mu(dw)$  is an almost periodic function on the strip $\{z\in\C:\,\a'<y<\b'\}$.
\smallskip

In our investigation we consider only almost periodic measures with supports in  closed horizontal strips of finite width. In this case the definition of almost periodic measures  can be simplified:
\medskip

{\bf Definition 4}(\cite{R}).
  A  measure $\mu$ with support in the closed strip $S_{[\a,\b]}=\{z=x+iy:-\infty<\a\le y\le\b<\infty\}$  is \emph{almost periodic}
  if for any continuous function $\p(z)$ with compact support in $\C$
 the convolution  $\p\star\mu(z)=\int\p(z-w)\mu(dw)$  is an almost periodic function on $\C$.
\medskip

{\bf Definition 5}.
 A measure $\mu$ on a closed horizontal strip $S_{[\a,\b]}$ is \emph{translation bounded}  if
$$
\sup_{t\in\R}|\mu|\{z\in S_{[\a,\b]}:\,t<\Re z<t+1\}<\infty.
$$

Clearly, a locally finite multiset $A=\{a_n\}$ is translation bounded iff the measure $\mu_A=\sum_n\d_{a_n}$ is translation bounded.
\begin{Pro}\label{P3}
Every almost periodic complex measure $\mu$  on $S_{[a,b]}$ is translation bounded.
\end{Pro}
An analog of this Proposition for the much wider class of almost periodic distributions in tube domains in $\C^n$ was proved in \cite{R}.  In our article we are interested only in
   complex measures on a closed strip of bounded width. In this case, there is a very short proof of Proposition \ref{P3}:
 \smallskip

{\bf Proof}. Let $X$ be the Banach space of all continuous functions $\p$ on $[0,1]\times[a,b]$ such that $\p(0+iy)=\p(1+iy)=0$ for all $y\in[a,\,b]$, let $\p$ be an arbitrary function from $X$,
and let $\tilde\p$ be a compactly supported continuous extension $\p$ on $\C$ such that $\tilde\p(x+iy)=0$ for $y\not\in[a-1,b+1]$. Since the convolution $\tilde\p\star\mu(z)$ is almost periodic in $\C$, we get that the acting
$$
 (\mu_t,\p)=\int_{\z\in S_{[a,b]}}\p(\z)\mu_t(d\z)=\int_{\z\in\C}\tilde\p(\z)\mu_t(d\z)=(\mu\star\tilde\p)(t)\quad\text{with}\quad \mu_t(\cdot)=\mu(t-\cdot)
$$
is  uniformly bounded in $t\in\R$. Applying  Banach--Steinhaus Theorem, we obtain the assertion. \qquad$\square$
\medskip

{\bf Definition 6} (see \cite{FRR}).
 A locally finite set $A=\{a_n\}\subset S_{[\a,\b]}$ (or $\subset S_{(\a,\b)}$) is \emph{almost periodic} if the measure $\mu_A=\sum_n\d_{a_n}$ is almost periodic in the same strip.
For $\a=\b=0$ we obtain the definition of almost periodic sets in $\R$.
\medskip

The original  definition,  due to M.Krein and B.Levin \cite[App. VI]{L}, looks as follows:
\medskip

{\bf Definition 7}.
A  locally finite set $A=\{a_n\}_{n\in\Z}\subset S_{[\a,\b]}$ is \emph{almost periodic} if for any $\e>0$ the set of its $\e$-almost periods
$$
E_\e=\{\tau\in\R:\,\exists\ \text{ a bijection } \s:\Z\to\Z \quad\text{such that}\ \sup_n |a_n+\tau-a_{\s(n)}|<\e\}
$$
has a nonempty intersection with every interval $(x,x+L_\e)$.
\medskip

The generalization of this definition to almost periodic sets in an open strip, in particular in $\C$, is due to H.Tornehave  \cite{T}.
In \cite{FRR} it was proved that Definition 7  and Tornhave's one are equivalent to Definition 6.

Note that the zero set of any holomorphic almost periodic function
on an open  strip is an almost periodic set in this strip. The converse is not true; the connection between almost periodic sets and zeros
of holomorphic almost periodic functions is rather complicated. A complete  description is given in \cite{F0} in terms of Chern cohomologies.
But in our article we will deal only with sets that {\it a priori} lie in a closed horizontal strip of finite width.
In this case every almost periodic set  is the zero set of an entire almost periodic function \cite{FRR}.

The zero set of Dirichlet series \eqref{S} is almost periodic and
lies in a horizontal strip of  finite width if and only if $\sup$ and $\inf$ of spectrum $Q$  belong to this spectrum (see \cite{L}, Ch.VI, Con.2).
In particular, the zero set of each exponential polynomial \eqref{s} lies in such a strip.
\medskip

Clearly, for every almost periodic set $A=\{a_n\}_{n\in\Z}\subset\SH$ the set $A'=\{\Re a_n\}$  is  almost periodic in $\R$. It was proved in Theorem 1 from \cite{F3} that
under condition $\Re a_n\le\Re a_{n+1}$ for all $n$
$$
\Re a_n=\r n+\psi(n)\quad\text{with an almost periodic mapping}\quad \psi:\Z\to\C\quad\text{and some density}\,1/\r.
$$
Therefore, with the appropriate numbering we obtain
\begin{Pro}\label{P4}
 For every almost periodic set $A=\{a_n\}_{n\in\Z}\subset\SH$  there is a bounded mapping $\phi:\,Z\to\C$ such that
$$
a_n=\r n+\phi(n).
$$
\end{Pro}

\bigskip
\section{Entire functions with almost periodic zeros}\label{S5}
\bigskip

By Hadamard's Theorem, every  entire function $g(z),\,z\in\C$, of exponential growth (i.e., $\log|g(z)|\le O(|z|)$ as $|z|\to\infty$) with zeros $a_n\in\C\setminus\{0\}$
has the form
\begin{equation}\label{g}
	g(z)= e^{dz}\prod_n(1-z/a_n)e^{z/a_n},\quad d\in\C.
\end{equation}
Note that each $a\in\C$  can be repeated  any finite number of times in the sequence $\{a_n\}$.

Further, zeros of $g$ satisfy the conditions
\begin{equation}\label{Lin1}
	\#\{n:\,|a_n|\le r\}=O(r),\qquad r\to\infty,
\end{equation}
and
\begin{equation}\label{Lin2}
	\sum_{n:\,|a_n|\le r}\frac{1}{a_n}=O(1),\qquad r\to\infty.
\end{equation}
On the other hand, if a sequence $A=\{a_n\}_{n\in\Z}\subset\C\setminus\{0\}$ satisfies these conditions,
then the function \eqref{g} is an entire function of exponential growth (Lindel\"of's Theorem  (see \cite{K})).

\begin{Pro}\label{P5}
Let $A=\{a_n\}_{n\in\Z}\subset\SH\setminus\{0\}$ be an almost periodic set. Then the points $a_n$ satisfy conditions \eqref{Lin1},  \eqref{Lin2}, and under suitable numeration the infinite product
\begin{equation}\label{f}
  f(z)=(1-z/a_0)\prod_{n\in\N} (1-z/a_n)(1-z/a_{-n})
\end{equation}
converges uniformly on compact sets in $\C$. Moreover, the sum
\begin{multline}\label{df}
\frac{f'(z)}{f(z)}=\frac{1}{z-a_0}+ \sum_{n\in\N}\left[\frac{1}{z-a_n}+
\frac{1}{z-a_{-n}}\right]=\\ \frac{1}{x+iy-\phi(0)}+\sum_{n\in\N}\left[\frac{1}{x+iy-\r n-\phi(n)}+\frac{1}{x+iy+\r n-\phi(-n)}\right]
\end{multline}
converges absolutely and uniformly (after discarding a finite number of members) on the sets $D_R:=\{z=x+iy:\,|x|\le R,\,|y|\ge M+1\}$, where $M:=1+\sup_{n\in\Z}|\phi(n)|$ and  $R<\infty$ is arbitrary.
\end{Pro}

{\bf Proof}. By Proposition \ref{P4}, we  can renumber $a_n$ such that $a_n=\r n+\phi(n)$, $n\in\Z$,
 with a bounded complex-valued function $\phi(n)$. Therefore the sum
$$
 \sum_{n\in\N}\left[\frac{1}{a_n}+\frac{1}{a_{-n}}\right]=\sum_{n\in\N}\left[\frac{1}{\r n+\phi(n)}+\frac{1}{-\r n+\phi(-n)}\right]
 $$
 converges absolutely, and the sum
 $$
 \sum_{n\in\N}\left[\log(1-z/a_n)+\log(1-z/a_{-n})\right]
$$
converges  absolutely and uniformly  on the sets $D_R$ after discarding a finite number of terms. The sums
$$
 \sum_{n\in\Z,|a_n|<r}\frac{1}{a_n} \quad\mbox{and}\quad \sum_{n\in\Z,0<|\r n|<r}\frac{1}{a_n}
 $$
 agree up to a uniformly bounded (with respect to $r$) number of terms, and each of these terms tends to $0$ as $r\to\infty$. Therefore the first sum has a finite limit as $r\to\infty$, and we obtain representation
 \eqref{f} and condition \eqref{Lin2}. Condition \eqref{Lin1} follows from boundedness of $\phi$.

 Since the sums
$$
\sum_{|\r n|<|x|+M+1}\left|\frac{1}{x+iy-\r n-\phi(n)}\right|+\sum_{\r n\ge|x|+M+1}\left|\frac{1}{x+iy-\r n-\phi(n)}+\frac{1}{x+iy+\r n-\phi(-n)}\right|
$$
are uniformly bounded for $z\in D_R$, we obtain \eqref{df}. $\square$

\begin{Pro}\label{P6}
Let $A$ be an almost periodic set in $\SH\setminus\{0\}$, and  the distribution $\hat\mu_A^c$ defined in \eqref{h} be a pure point measure  \eqref{b} such that
\begin{equation}\label{var0}
  |\hat\mu_A^c|(-r,r)=\sum_{|\g|<r}|b_\g|=O\left(e^{Lr}\right),\qquad r\to\infty,
\end{equation}
with  $L>2\pi H$. Then the function \eqref{f} satisfies the condition
\begin{equation}\label{der}
f'(\z)/f(\z)=-2\pi i\sum_{\g\in\G\cap(0,+\infty)}b_\g e^{2\pi i\g\z}-\pi ib_0,
\end{equation}
and the sum in the right-hand side absolutely converges on every horizontal line with $\Im\z>L/2\pi$.
\end{Pro}
{\bf Proof}.
Fix $\z\in\C,\,\Im\z=K/2\pi$ with $K>L$.  Set $e_\z(t)=-2\pi ie^{2\pi it\z}$ for $t>0$ and $e_\z(t)=0$ for $t\le0$.
Clearly, its c-Fourier transform
$$
\hat e^c_\z(z)=-\int_0^\infty 2\pi i e^{2\pi it\z}e^{-2\pi izt}dt=1/(\z-z)
$$
is well-defined for $z\in\SH$.
Let $\p(t)$ be an even nonnegative $C^\infty$-function such that $\supp\p\subset(-1,1)$ and $\int\p(t)dt=1$. Set  $\p_\e(t)=\e^{-1}\p(t/\e)$ for $\e>0$. We get
$$
 \widehat{(\p_\e(z))^c}=\hat\p^c(\e z),\quad |\hat\p^c(\e z)|\le e^{2\pi H},\quad\hat\p^c(\e z)\to1\,\text{ as }\,\e\to0.
$$
The functions $e_\z(t),\,e_\z\star\p_\e(t)$ and all their derivatives satisfy \eqref{e}. By Proposition \ref{P3} we have $\mu_A(-r,r)=O(r)$.
It follows from Proposition \ref{P1}  that \eqref{h} holds for the function $e_\z\star\p_\e$. The equality \eqref{com} yields
\begin{equation}\label{d}
   (\hat\mu_A^c,e_\z\star\p_\e)=(\mu_A,\hat e^c_\z\hat\p_\e^c)
\end{equation}
We have
\begin{equation}\label{c}
 (\mu_A(z),\hat e^c_\z(z)\hat\p_\e^c(z))=\frac{\hat\p_\e^c(a_0)}{\z-a_0}+\sum_{n\in\N}\hat\p_\e^c(a_{-n})\left[\frac{1}{\z-a_n}+\frac{1}{\z-a_{-n}}\right]
 +\sum_{n\in\N}\frac{\hat\p^c(\e a_n)-\hat\p^c(\e a_{-n})}{\z-a_n}.
\end{equation}
 By Proposition \ref{P5},  the first sum is bounded uniformly in $\e$. Then,
$$
   |\z-a_n|\ge\max\{|\Im \z-\Im\phi(n)|,|\Re \z-\r n-\Re \phi(n)|\}\ge\max\{K/2\pi-H,||\Re\z-\r n|-M|\}.
$$
The  series
$$
   \sum_{n\in\Z}[\max\{K/2\pi-H,||\Re\z-\r n|-M|\}]^{-2}
$$
 represents a periodic function on the line $\Im\z=K/2\pi$, hence,
$$
   S:=\sup\left\{\sum_{n\in\N}|\z-a_n|^{-2}:\,\Im\z=K/2\pi\right\}<\infty.
$$
Using Cauchy-Schwarz-Bunyakovskii inequality, we obtain
 \begin{equation}\label{E}
 \left|\sum_{n\in\N}\frac{\hat\p^c(\e a_n)-\hat\p^c(\e a_{-n})}{\z-a_n}\right|\le S^{1/2}\left(\sum_{n\in\N}|\hat\p^c(\e a_n)-\hat\p^c(\e a_{-n})|^2\right)^{1/2}
 \end{equation}
 Next, taking into account that $\hat\p^c$ is even,  we get
$$
\hat\p^c(\e a_n)-\hat\p^c(\e a_{-n})\le\hat\p^c(\e\r n+\e\phi(n))-\hat\p^c(\e\r n-\e\phi(-n))\le2M\e\max\{|(\hat\p^c)'(z)|:\,z\in I_n\},
$$
where  the number $M$ is the same as in Proposition \ref{P5} and $I_n$ is the segment $[\e\r n-\e\phi(-n),\e\r n+\e\phi(n)]$.
It follows from  \eqref{N} that
$$
\max\{|(\hat\p^c)'(z)|:\,z\in I_n\}\le C\max_{I_n}\min\{1,|z|^{-1}\}\le\min\{1,(\e\r n-\e M)^{-1}\}.
$$
 Hence,
 $$
 \sum_{n\in\N}|\hat\p^c(\e a_n)-\hat\p^c(\e a_{-n})|^2\le  4\e^2M^2C^2\left[ \sum_{1\le n\le 1/\e\r+M/\r}1+\sum_{n>1/\e\r+M/\r}(\e\r n-\e M)^{-2}\right].
 $$
The first sum in the right-hand side is $O(1/\e)$, and second one is $o(1/\e^2)$. Therefore, \eqref{E} tends to zero as $\e\to0$, and \eqref{c} is uniformly bounded for small $\e$.
Applying Lebesgue's Dominate Convergence Theorem, we can take the limit in \eqref{c} as $\e\to0$. We obtain
\begin{equation}\label{q1}
(\mu_A(z),\hat e^c_\z(z)\hat\p_\e^c(z))\to \frac{1}{\z-a_0}+\sum_{n\in\N}\left[\frac{1}{\z-a_n}+\frac{1}{\z-a_{-n}}\right],\qquad\e\to0.
\end{equation}

 Set
$$
n(s)=\sum_{\g\in\G:\,0<\g<s}|b_\g|.
$$
By \eqref{var0},  $n(s)\le Ce^{Ls}$. We have for $y>L/2\pi$
\begin{equation}\label{i1}
   \sum_{\g\ge r}|b_\g|e^{-2\pi\g y}=\int_r^\infty e^{-2\pi sy}n(ds)\le\lim_{T\to\infty}n(T)e^{-2\pi Ty}+2\pi y\int_r^\infty e^{-2\pi sy}n(s)ds<\infty.
\end{equation}
Therefore the series $\sum_{\g>0}b_\g e^{2\pi i\g\z}$ converges absolutely and uniformly on the line $\Im\z=K/2\pi$.

Next, we have for $\e>0$
\begin{equation}\label{i3}
 \frac{i}{2\pi}(\hat\mu_A^c(t),e_\z\star\p_\e(t))= \sum_{\g\in\G}b_\g e^{2\pi i\g\z}\int_{-\e}^\g e^{-2\pi is\z}\p_\e(s)ds.
 \end{equation}
It is easy to see that all  integral on the right-hand side do not  exceed $2e^{\e K}\max\p$, and for $\e\to0$
$$
\int_{-\e}^\g e^{-2\pi is\z}\p_\e(s)ds=\int_{-1}^{\min\{\g/\e,1\}}e^{-2\pi i\e s\z}\p(s)ds\to\begin{cases}1,&\text{if }\g>0,\\1/2,&\text{if }\g=0,\end{cases}                                                                                  $$
$$
\int_{-\e}^\g e^{-2\pi is\z}\p_\e(s)ds=\int_{-1}^{\max\{\g/\e,-1\}}e^{-2\pi i\e s\z}\p(s)ds\to0,\quad\text{if }\g<0.
$$
Using Lebesgue's Dominate Convergence Theorem, we obtain from \eqref{i3}
$$
  (\hat\mu_A^c(t),e_\z\star\p_\e(t))\to-2\pi i\sum_{\g\in\G\cap(0,+\infty)}b_\g e^{2\pi i\g\z}-\pi ib_0,\qquad\e\to0.
$$
 \eqref{df},   \eqref{d},  and \eqref{q1} yield \eqref{der}. $\square$

\bigskip
 \section{Proof of  Theorem \ref{T1}}\label{S6}
\bigskip

Since the zero set $A$ of the Dirichlet series \eqref{S} lies in a strip of finite width, we get $\sup\O\in\O$ and $\inf\O\in\O$ (see Section \ref{S4}).
Without loss of generality suppose $\sup\O=\k,\,\inf\O=-\k$. We have $q_\k\neq0,\,q_{-\k}\neq0$. Then
$$
   Q(z)=q_{-\k} e^{-2\pi i\k z}(1+P(z)),\qquad P(z)=\sum_{\o\in\O\setminus\{-\k\}}\frac{q_\o}{q_{-\k}e^{2\pi(\k+\o)y}}e^{2\pi i(\k+\o)x},\quad z=x+iy.
$$
 Taking into account that $\sum_\o|q_\o|<\infty$, we choose a finite number of elements $\o_1,\dots,\o_N\in\O\setminus\{-\k\}$ and then $s^*>H$ such that
$$
\sum_{\o\in\O\setminus\{-\k,\o_1,\dots,\o_N\}}|q_\o/q_{-\k}|<1/3,\qquad \sum_{j=1}^N e^{-2\pi(\o_j+\k)s^*} |q_{\o_j}/q_{-\k}|<1/3.
$$
So $\|P(x+is^*)\|_\W<2/3$, and by Proposition \ref{P2},
\begin{equation}\label{sn1}
\log(1+P(x+is^*))=\sum_{\g\in\G^*}p_\g e^{2\pi i\g x}\quad\text{with}\quad\G^*\subset(0,+\infty),\quad\sum_{o\in\G^*}|p_\g|<\infty.
\end{equation}
Therefore,
\begin{equation}\label{sn2}
Q'(x+is^*)/Q(x+is^*)=[\log Q(x+is^*)]'=\sum_{\g\in\G^*}2\pi i\g p_\g e^{2\pi i\g x}-2\pi i\k.
\end{equation}
Since $Q$ is an almost periodic function, we can apply  Lemma 1  \cite[Ch.6]{L}, and for any $\e>0$ and $s<\infty$ find a number $m=m(\e,s)>0$ such that
\begin{equation}\label{inf}
 |Q(z)|\ge m\quad\mbox{for}\quad |\Im z|\le s\quad\mbox{and}\quad \dist(z,A)\ge\e.
\end{equation}
Hence,
$$
\inf_{x\in\R}|Q(x\pm is)|>0
$$
for every fixed $s>H$, and $1/Q(x+is^*)\in\W$. Since $\O$ is bounded, we get $Q'(x+is^*)\in\W$ and $(Q'/Q)(x+is^*)\in\W$. Therefore,
\begin{equation}\label{sn3}
\sum_{o\in\G^*}|\g p_\g|<\infty.
\end{equation}
The same arguments show that for some countable set $\G_*\subset(-\infty,0)$ and $s_*$ large enough
\begin{equation}\label{sn4}
 Q'(x-is_*)/Q(x-is_*)=\sum_{\g\in\G_*}2\pi i\g p_\g e^{2\pi i\g x}+2\pi i\k,\quad \sum_{\g\in\G_*}| p_\g|<\infty,\quad \sum_{\g\in\G_*}|\g p_\g|<\infty.
\end{equation}
In what follows we set $s=\max\{s^*,s_*\}$.
\medskip

Since $Q$ is an almost periodic function, we see that its zero set $A\subset\SH$  is almost periodic.  Therefore by Proposition \ref{P3},  the numbers $\#\{a_n\in A:\,x<\Re a_n<x+1\}$ are bounded uniformly in $x\in\R$.
Consequently, for $\e$ small enough every interval $(x,x+1)$ contains a number $L$ such that $\dist(z,A)\ge\e$ for all points of the segment $[L-is,L+is]$.  By \eqref{inf}, there exist $m=m(s,\e)>0$ and two sequences
$L_k\to+\infty,\,L'_k\to-\infty$ such that
$$
   |Q(x+iy)|\ge m>0\quad\text{for}\quad x=L_k \quad\text{or}\quad x=L'_k,\quad |y|\le s.
$$
Set $\p\in\D$.  The function $\hat\p^c(z)$  is holomorphic in $\C$ and by \eqref{N}, the integrals of the function $\hat\p^c(z)Q'(z)Q^{-1}(z)$ over  boundaries of the rectangles $\{z:\,L'_k<x<L_k,-s<y<s\}$ tend to
$$
   \int_{-\infty}^{+\infty}\hat\p^c(x-is)Q'(x-is)Q^{-1}(x-is)dx-\int_{-\infty}^{+\infty}\hat\p^c(x+is)Q'(x+is)Q^{-1}(x+is)dx=:I^--I^+.
$$
 Using the Residue Theorem, we get
 \begin{equation}\label{r}
\frac{I^--I^+}{2\pi i}=\sum_{\l:Q(\l)=0}\Res_\l \hat\p^c(z)\frac{Q'(z)}{Q(z)}
 =\sum_{\l:Q(\l)=0}a(\l)\hat\p^c(\l)= (\mu_A,\hat\p^c),
\end{equation}
where $a(\l)$ is the multiplicity of the zero at the point $\l$.
  Using \eqref{sn2}, \eqref{sn4}, we get
\begin{multline*}
 \frac{I^--I^+}{2\pi i}=\sum_{\g\in\G_*}\g p_\g\int_{-\infty}^{+\infty}\hat\p^c(x-is)e^{2\pi i\g x}dx-\sum_{\g\in\G^*}\g p_\g\int_{-\infty}^{+\infty}\hat\p^c(x+is)e^{2\pi i\g x}dx\\
+\k\int_{-\infty}^{+\infty}\hat\p^c(x-is)dx+\k\int_{-\infty}^{+\infty}\hat\p^c(x+is)dx.
\end{multline*}
To calculate the inverse Fourier transform, we obtain. by \eqref{ext},
 \begin{equation}\label{i}
  \frac{I^--I^+}{2\pi i}=\sum_{\g\in\G_*}e^{-2\pi\g s}\g p_\g \p(\g)-\sum_{\g\in\G^*}e^{2\pi\g s}\g p_\g\p(\g)+2\k\p(0).
 \end{equation}
 Set $\G=\G_*\cup\G^*\cup\{0\},\quad b_\g=-\g p_\g e^{2\pi\g s}\quad\text{for}\ \g\in\G^*,\quad b_\g=\g p_\g e^{-2\pi\g s}\quad\text{for}\ \g\in\G_*,\quad b_0=2\k$.
\medskip

It follows from \eqref{r} and \eqref{i} that $\nu=\sum_{\g\in\G} b_\g\d_\g$   is a measure on $\R$ such that
$$
   (\nu,\p)=\sum_{\g\in\G}b_\g\p(\g)=(\mu_A,\hat\p^c).
$$
 Since this equality  is valid for every $\p\in\D$, we obtain $\hat\mu_A^c=\nu$.

If $\supp\p\subset(-r,r)$ for some $r<\infty$, then  from the definitions $\nu$ and $b_\g$, we obtain
$$
|(\hat\mu_A^c,\p)|\le\sum_{\g\in\G\cap[-r,r]}|b_\g|\le e^{2\pi rs}\left(\sum_{\g\in\G_*}|\g p_\g|+\sum_{\g\in\G^*}|\g p_\g|+2\k\right)\sup_{|t|\le r}|\p(t)|.
$$
Taking into account \eqref{sn3} and \eqref{sn4}, we obtain the estimate
$$
  |(\hat\mu_A^c,\p)|\le  C(s)e^{2\pi sr}\sup_{|t|\le r}|\p(t)|,\quad C(s)<\infty.
$$
This bound  remains true for all continuous functions on $[-r,r]$ that vanish at $\pm r$. Therefore,
$$
  |\hat\mu_A^c|(-r,r)\le C(s)e^{2\pi sr}.
$$
Property \eqref{var} is proved. Property \eqref{neig} follows from \eqref{sn1},  \eqref{sn4}, and the definition of $b_\g$.

If $Q$ is an exponential polynomial, then it follows from Proposition \ref{P2}  that  $\G^*$ and $\G_*$ are locally finite.\  $\square$

{\bf Remark}. Without loss of generality assume $0\not\in A$. Then the function \eqref{f} is an entire function of exponential growth with the zero set $A$, and $Q(z)$ is the same.
Consequently, $Q(z)=Ce^{dz}f(z),\,d,C\in\C$. This equality and Proposition \ref{P6} yield
\begin{multline*}
  \frac{Q'(x+is)}{Q(x+is)}+2\pi i\k=\frac{f'(x+is)}{f(x+is)}+d+2\pi i\k=-2\pi i\sum_{\g\in\G\cap\R_+}\frac{b_\g}{e^{2\pi\g s}} e^{2\pi i\g x}+i(2\pi\k-\pi b_0)+d.
\end{multline*}
Since $\spec\{Q'/Q+2\pi i\k\}\subset(0,+\infty)$, we obtain  $\Re d=0$, $\Im d=\pi b_0-2\pi\k$, and $Q(z)=Cf(z)e^{(\pi b_0-2\pi\k)iz}$.

\bigskip
 \section{Proof of  Theorem \ref{T2}}\label{S7}
\bigskip

Let $A\subset\SH\setminus\{0\}$ be an almost periodic set, $\mu_A$ be the measure \eqref{a0}, and the pure point measure $\hat\mu_A^c$ satisfy \eqref{neig} and \eqref{var0} with $L>2\pi H$.
Without loss of generality, assume that $0\not\in A$. We integrate the equality \eqref{der} over the segment $[iy_0,x+iy_0]$ with $y_0>L/2\pi$
 and change the order of summation and integration.  We get
$$
  \log f(x+iy_0)-\log f(iy_0)=\int_{iy_0}^{x+iy_0}\frac{f'(\zeta)}{f(\zeta)}d\zeta=-\sum_{\g\in\G\cap(0,+\infty)}b_\g \frac{(e^{2\pi i\g x}-1)e^{-2\pi\g y_0}}{\g}-ib_0\pi x.
$$
Hence,
$$
\log f(x+iy_0)+ib_0\pi x= -\sum_{\g\in\G\cap(0,+\infty)}\frac{b_\g}{\g e^{2\pi\g y_0}}e^{2\pi i\g x}+\sum_{\g\in\G\cap(0,+\infty)}\frac{b_\g}{\g e^{2\pi\g y_0}}+\log|f(iy_0)|.
$$
It follows from \eqref{neig} and \eqref{i1} that the  sums
$$
  \sum_{0<\g<1}|b_\g|\g^{-1}e^{-2\pi\g y_0}\quad\text{and}\quad \sum_{\g\ge1}|b_\g|\g^{-1}e^{-2\pi\g y_0}
$$
converge, hence the function
 $$
   f(x+iy_0)e^{ib_0\pi x}=\exp\{\log f(x+iy_0)+ib_0\pi x\}
 $$
 belongs to $\W$. Therefore,
 $$
   f(x+iy_0)e^{ib_0\pi x}=\sum_{\o\in\O}\b_\o e^{2\pi i\o x},\qquad \sum_{\o\in\O}|\b_\o|<\infty,
 $$
with $\b_\o\in\C$ and a countable spectrum $\O$.  The entire function $g(z)=f(z+iy_0)e^{ib_0\pi z}$ has the exponential growth and is bounded on the line $z=x\in\R$.
By the Phragmen--Lindelof Principle, it is bounded on every horizontal strip of a finite width. It follows from Section \ref{S4} that $g(z)$ is almost periodic function on $S_{(\infty,+\infty)}$.
By \cite[\S 1, Ch.VI]{L}, it follows that $\O$ is bounded. Hence the function
$$
    f(z)=\sum_{\o\in\O}\b_\o e^{\pi(2\o-b_0)y_0}e^{\pi i(2\o-b_0)z}
$$
is  Dirichlet series of the form \eqref{S}.

\medskip
If $\G$ is locally finite, then the function $\log f(x+iy_0)+ib_0\pi x\in\W$ has nonnegative locally finite spectrum. By Proposition \ref{P2},
the same is valid for the function $f(x+iy_0)e^{ib_0\pi x}$. Since $\O$ is bounded, we see that only a finite number of coefficients $\b_\o$ does not vanish.
Hence, $f(z)$ is an exponential polynomial. \  $\square$
\medskip

I would like to thank the reviewer for numerous comments that allowed me to improve the article. I am also grateful to Professor Szilard Revesz from Renyi Institute of Mathematics
for his hospitality, for attention to my research, and for useful discussions.

\end{document}